\input amstex
\documentstyle{amsppt}
\loadbold
\pageheight{7.0in}

\magnification=\magstep 1
\CenteredTagsOnSplits
\NoBlackBoxes
\widestnumber \key{Pr98a}
%
\def\curraddr#1\endcurraddr{\address {\it Current address\/}: #1\endaddress}

\catcode`\@=11
\redefine\logo@{}
\catcode`\@=13

\def\<{\left<}
\def\>{\right>}

\topmatter
\title Cocycle conjugacy classes of binary shifts
\endtitle
\author
Geoffrey L. Price
\endauthor
\affil
United States Naval Academy
\endaffil

\address
Department of Mathematics,
United States Naval Academy,
Annapolis MD 21402, U.S.A.
\endaddress

\email
glp\@usna.edu
\endemail

\abstract We show that every binary shift on the hyperfinite $II_1$ factor $R$ is cocycle conjugate to at least countably many non-conjugate binary shifts.  This holds in particular for binary shifts of infinite commutant index.
\endabstract

\date
February 12, 2016
\enddate

\dedicatory
In memory of William B. Arveson
\enddedicatory

\endtopmatter

\document

\subhead
1. Preliminaries
\endsubhead

Let $(a) = a_0,a_1,a_2,\dots$ be a {\it bitstream}, i.e. a sequence of $0$'s and $1$'s in $GF(2)$.  Let $v_0,v_1,v_2,\dots$ be a sequence of self-adjoint unitary operators (which we shall call generators) that satisfy the translation-invariant commutation relations

$$
v_iv_{i+j} = (-1)^{a_j}v_{i+j}v_{i}
$$
for all non-negative integers $i$ and $j$.  The commutation relations imply that $a_0 = 0$.  The norm closure of the set of all linear combinations of the words in the $v_i$'s is isomorphic to the CAR algebra, $\frak{A}$, if and only if the mirror sequence $\dots , a_2,a_1,a_0,a_1,a_2,\dots$ is not periodic, \cite{Pr87, Theorem 2.3}.  Using the GNS representation for the unique trace on $\frak{A}$, \cite{KR}, which sends all non-trivial words in the $v_i$'s to $0$, the hyperfinite $II_1$ factor $R$ is obtained as the strong closure of the set of all linear combinations of the words in the $v_i$'s.  Then the mapping sending $v_i$ to $v_{i+1}$ gives rise to a unital $*$-homomorphism $\alpha$ on $R$, called a {\it binary shift}.  The image $\alpha(R)$ has subfactor index $2$, i.e. $[R: \alpha(R)]=2$.  Note that $R$ is generated by $\alpha(R)$ and $v_0$.

\definition{Definition 1, cf. \cite{Po}}  Two unital $*$-endomorphisms $\alpha$ and $\beta$ on $R$ are said to be {\it conjugate} if there is a $*$-automorphism $\gamma$ of $R$ such that
$$
\alpha(v) = \gamma(\beta(\gamma^{-1}(v))), \text{for all}\,\, v\in R.
$$
Two unital $*$-endomorphisms $\alpha$ and $\beta$ are said to be {\it cocycle conjugate} if there is a unitary $u\in R$ so that $Ad(u)\circ \alpha$ and $\beta$ are conjugate.
\enddefinition

\proclaim{Theorem 1} \cite {Po, Theorem 3.6} Binary shifts $\alpha$ and $\beta$ on $R$ are conjugate if and only if they have the same bitstream.
\endproclaim

For $n\in \Bbb N$ let $\Cal A_n$ be the $n\times n$ Toeplitz matrix over $GF(2)$, given by

$$
\Cal A_n  = \bmatrix a_0 & a_1 & a_2 & a_3 & \hdots & a_{n-1}
\\
a_1 & a_0 & a_1 & a_2 & \hdots & a_{n-2}
\\
a_2 & a_1 & a_0 & a_1 & \hdots & a_{n-3}
\\
\vdots & \vdots & \vdots & \vdots & \ddots & \vdots
\\
a_{n-1} & a_{n-2} & a_{n-3} & \hdots & \hdots & a_0
\endbmatrix
$$

Note that $\Cal A_n$ is symmetric with main diagonal all $0$'s.  For each $n\in \Bbb N$ let $\nu_{n} = \text{null}(\Cal A_n)$, i.e. the dimension of $ker(\Cal A_n)$.  The importance of the nullity sequence $\{\nu_n: n\in \Bbb N\}$ is that it gives information about the centers of the algebras $\frak{A}_n$ generated by $v_0,\dots,v_{n-1}$ (see Theorem $3$ and Theorem $4$ below).

\proclaim{Theorem 2} \cite{PP, Theorem 5.4},\cite{CP, Theorem 2.7}  The sequence $\{\nu_n: n\in \Bbb N \}$ is the concatenation of finite length strings of the form $1,0$ or $1,2,\dots,r-1,r,r-1,\dots ,0$ for some $r > 1$, where the length $2r$ of the string may vary.
\endproclaim

\remark{Remark 1}  If the mirror bitstream is periodic there is a non-negative integer $n$ such that $\nu_{n+j+1} = \nu_{n+j}+1$ for all $j\in \Bbb N$,\cite{Pr98}, and in this case the $C^*$-algebra $C^*(v_0,v_1,\dots)$ is isomorphic to $\Cal T \otimes C(X)$, where $X$ is the Cantor set, and $\Cal T$ is either $\Bbb CI$ or the tensor product of finitely many copies of $M_2(\Bbb C)$, \cite{AP, Theorem 3.1}.
\endremark

\proclaim{Theorem 3}  \cite{PP, Corollary 5.5},\cite{Pr98, Lemma 3.2}  Let $\alpha$ be a binary shift on $R$ with bitstream $(a)$.  For any positive integer $n$ there exists a one-to-one correspondence between vectors $\bold c= [c_0,c_1,\dots,c_{n-1}]^t$  in $ker(\Cal A_n)$ and ordered words $v_0^{c_0}v_1^{c_1}\cdots v_{n-1}^{c_{n-1}}$ in the center $\frak{Z}_n$ of $\frak{A}_n$.  Therefore, as $\frak{Z}_n$ consists of linear combinations of the words it contains, its dimension as an algebra over $\Bbb C$ is $2^{\nu_n}$.  In particular, if $\nu_n =0$ then $\frak{Z}_n = \Bbb CI$ and $\frak{A}_n$ is a full matrix algebra over $\Bbb C$.
\endproclaim

The following result about the centers $\frak{Z}_n$ of the algebras $\frak{A}_n$, along with the correspondence between $\frak{Z}_n$ and $ker(\Cal A_n)$ discussed above, shows that to understand the structure of $ker(\Cal A_n)$ for all $n$ it suffices to know the structure of $ker(\Cal A_n)$, where $\nu_{n-1} =0$ and, (therefore necessarily, by Theorem $2$), $\nu_{n} =1$.

\proclaim{Theorem 4} \cite{PP, Lemma 6.5},\cite{Pr98, Theorem 3.4}  Let $n\in \Bbb N$ and $d\in \Bbb N$ be such that the nullities $\nu_n$ through $\nu_{n+d-1}$ are $1$ through $d$, respectively, and $\nu_{n+d}$ through $\nu_{n+2d-1}$ are $d-1$ through $0$, respectively.  Let $z=v_0^{c_0}v_1^{c_1}\cdots v_{n-1}^{c_{n-1}}$ be the word generating $\frak{Z}_n$, then the exponents $c_0,c_1,\dots,c_{n-1}$ form a palindrome, i.e. $[c_0,c_1,\dots,c_{n-1}] = [c_{n-1},c_{n-2},\dots, c_0]$, with $c_0 = c_{n-1}=1$.  For $m = 0, 1,\dots,d-1$, the center
$\frak{Z}_{n+m}$ is generated by the words $z,\alpha(z),\cdots,\alpha^m(z)$.  $\frak{Z}_{n+d}$ is generated by $\alpha(z),\dots,\alpha^{d-1}(z)$,
$\frak{Z}_{n+d+1}$ by $\alpha^2(z),\dots,\alpha^{d-1}(z)$, and so on, and $\frak{Z}_{n+2d-1}$ is trivial.
\endproclaim
\subhead
2. Main Result
\endsubhead
We now construct countably many binary shifts that are cocycle conjugate to a given binary shift $\alpha$.  Fix $n\in \Bbb N$ such that $\nu_{n-1}=0$ and $\nu_n=1$.  As above, let $z=v_0^{c_0}v_1^{c_1}\cdots v_{n-1}^{c_{n-1}}$ be the word generating $\frak{Z}_n$, the center of the algebra $\frak{A}_n$ generated by $v_0,v_1,\dots,v_{n-1}$. By Theorem $4$, the exponents of $z$ form a palindrome.  Define a unitary operator $u$ by

$$
u = \left\{
\aligned
\frac{I+z}{\sqrt{2}},\,\,\, &\text{if}\,\,\, z^* = -z,\\
\frac{I+iz}{\sqrt{2}},\,\,\,&\text{if} \,\,\,z^* = z.
\endaligned
\right.
$$

Note for any word $w$ in the generators,

$$
Ad(u)w = \left\{
\aligned
w,\,\, &\text{if}\,\,\, zw = wz, \\
-zw, \,\,\,&\text{if}\,\,\, zw = -wz \,\,\,\text{and}\,\,\,z^* = -z,\\
-izw, \,\,\,&\text{if}\,\,\, zw = -wz \,\,\,\text{and}\,\,\,z^* = z.
\endaligned
\right.
$$

Set $\beta = Ad(u)\circ \alpha, u_0 = v_0$ and for $i\in \Bbb N,$ let $u_i = \beta^i(v_0)$.  Since $z$ commutes with $v_0,v_1,\dots,v_{n-1}$ it follows that $u_i = v_i$ for $i=0,1,\dots,n-1$.

First suppose $\nu_{n+1} = 0$.  Then $\frak{Z}_{n+1}$ is trivial.  Since $z\in \frak{Z}_{n}$, and therefore commutes with $v_0$ through $v_{n-1}$, we conclude that $z$ anticommutes with $v_{n}$ (otherwise $z$ would be in $\frak{Z}_{n+1}$, which is false).  So we have

$$
u_n = \beta^n(v_0) = \beta(\beta^{n-1}(v_0)) = \beta(v_{n-1}) = u^*v_nu = \lambda_0 zv_n, \tag1
$$
(where $\lambda_0$ is $-1$ or $-i$), and so $u_{n+1}=\beta(u_n)=u^*(\lambda_0\alpha(z)v_{n+1})u$ is either $\lambda_1 z\alpha(z)v_{n+1}$ or $\lambda_1 \alpha(z)v_{n+1}$ with $\lambda_1$ a fourth root of unity.  Write $u_{n+1}$ as $ \lambda_1 z^{r_0}\alpha(z)v_{n+1}$.
Using the symmetry arising from the fact that the exponent pattern of $z$ forms a palindrome, we conclude that since $z$ anticommutes with $v_n$, then
$\alpha(z)$ anticommutes with $v_0$.  Also $v_0$ commutes with $z$.  Since $u_0 = v_0$ and since $u_{n+1}=\lambda_1 z^{r_0}\alpha(z)v_{n+1}$, we conclude that $b_{n+1}\neq a_{n+1}$.  Therefore $\alpha$ and $\beta$ have different bitstreams.

Finally observe that the operators $u_0,u_1,\dots,u_{n-1},u_n,u_{n+1},u_{n+2},\dots$ are the operators $v_0,v_1,\dots,v_{n-1},\lambda_0zv_n,\lambda_1 z^{r_0}\alpha(z)v_{n+1},\lambda_2 z^{r_1}\alpha(z)^{r_0}\alpha^2(z)v_{n+2},\dots$ for some scalars $\lambda_i$ of modulus $1$.  Since $z\in \frak{A}_n$ it is clear that for all $m\in \Bbb N$, $\frak{B}_m = \frak{A}_m$, where $\frak{B}_m$ is the algebra generated by $u_0,\dots,u_{m-1}$.  So $\beta$ is a binary shift (with unitary generators $u_0,u_1,\dots$ and bitstream $(b) = b_0,b_1,\dots$) which is cocycle conjugate to $\alpha$.  Since $\alpha$ and $\beta$ have different bitstreams they are not conjugate, by Theorem $1$.

Next suppose the nullity sequence for $\alpha$ satisfies the hypotheses of Theorem $4$ with $d\geq 2$. Then by the conclusion of the theorem, $z$, the generator of $\frak{Z}_n$, is in $\frak{Z}_{n+d-1}$ but not $\frak{Z}_{n+d}$.   So $z$ commutes with generators $v_0$ through $v_{n+d-2}$ and anticommutes with $v_{n+d-1}$.  Setting $u$ as above and $\beta = Ad(u)\circ \alpha$ we see that $u_0$ through $u_{n+d-2}$ coincide with $v_0$ through $v_{n+d-2}$, but $u_{n+d-1}= \beta(u_{n+d-2}) = \beta(v_{n+d-2}) = u^*v_{n+d-1}u = \lambda_0zv_{n+d-1}$, $u_{n+d} = u^*\alpha(\lambda_0zv_{n+d-1})u =\lambda_1 z^{r_0}\alpha(z)v_{n+d}$, $u_{n+d+1} = \lambda_2 z^{r_1}\alpha(z)^{r_0}\alpha^2(z)v_{n+d+1}$, and so on, for scalars $\lambda_i$ of modulus $1$.

As $\alpha^{d-1}(z)$ is in $\frak{Z}_{n+2d-2}$, and therefore commutes with $v_0$ through $v_{n+2d-3}$, it follows that $\alpha^d(z)$ commutes with $v_1$ through $v_{n+2d-2}$.  As $\frak{Z}_{n+2d-1}$ is trivial, however, and as $\frak{A}_{n+2d-1}$ is generated by $v_0$ through $v_{n+2d-2}$ we conclude that $\alpha^d(z)$ anticommutes with $v_0$.  Therefore the following statements are true:

$$
\align
&u_0,u_1,\dots,u_{n+d-2}\,\, \text{agree with}\,\, v_0,v_1,\dots,v_{n+d-2},\,\, \text{respectively},\\
&u_{n+d-1}\,\,\text{is a scalar multiple of}\,\,zv_{n+d-1}\\
&u_{n+d}\,\,\text{is a scalar multiple of}\,\,z^{s_0}\alpha(z)v_{n+d},\,\, \text{for some exponent}\,\,s_0,\\
&u_{n+d+1}\,\,\text{is a scalar multiple of}\,\,z^{s_1}\alpha(z^{s_0})\alpha^2(z)v_{n+d+1},\,\, \text{for some exponent}\,\,s_1,\\
\vdots \\
&u_{n+2d-2}\,\,\text{is a scalar multiple of}\,\,z^{s_{d-2}}\cdots\alpha^{d-2}(z)^{s_0}\alpha^{d-1}(z)v_{n+2d-2},\,\, \text{and}\\
&u_{n+2d-1}\,\,\text{is a scalar multiple of}\,\,z^{s_{d-1}}\cdots\alpha^{d-1}(z)^{s_0}\alpha^{d}(z)v_{n+2d-1}.\\
\endalign
$$

It follows that the bitstream entries $b_0$ through $b_{n+2d-2}$ coincide with $a_0$ through $a_{n+2d-1}$ but, as
$u_{n+2d-1}$ is a scalar multiple of $z^{s_{d-1}}\cdots\alpha^{d-2}(z^{s_1})\alpha^{d-1}(z^{s_0})\alpha^d(z)v_{n+2d-1}$, and $u_0$ anticommutes with $\alpha^d(z)$, it follows that $b_{n+2d-1}\neq a_{n+2d-1}$.

An argument similar to the one given in the original case shows that $\frak{A}_m = \frak{B}_m$ for all $m$, so $\beta$ is a binary shift on $R$ cocycle conjugate, but not conjugate, to $\alpha$. So we have nearly completed the proof of the following.

\proclaim{Theorem 5}  Any binary shift on the hyperfinite $II_1$ factor $R$ is cocycle conjugate to at least countably many others.
\endproclaim
\demo{End of proof}  Let $\beta_1$ be the binary shift $\beta$ constructed above.  To construct $\beta_2$ from $\alpha$, choose $n_2\in \Bbb N$ such that $n_2 > n_1 = n$ and $\nu_{n_2-1} = 0, \nu_{n_2}=1$, and mimic the construction already made for $\beta$.  Note that the first place where the bitstreams for $\beta_2$ and $\alpha$ differ will occur past the first place where the bitstream for $\beta_1$ and $\alpha$ differ, so $\beta_1$ and $\beta_2$ are not conjugate. Continuing this process we can construct countably many binary shifts which are mutually non-conjugate but which are all cocycle conjugate to $\alpha$.
\enddemo

\definition{Definition 2}  For $k\in \Bbb N$, a binary shift $\alpha$ on $R$ is said to have commutant index $k$ if the relative commutant algebra $\alpha^k(R)^\prime \cap R$ is non-trivial and $k$ is the first non-negative integer for which this is the case.  If $\alpha^k(R)^\prime \cap R = \Bbb CI$ for all $k\in \Bbb N$ then we say that $\alpha$ has infinite commutant index.
\enddefinition

Here are a few facts about the commutant index of a binary shift.  The minimal commutant index is $2$, \cite{J}, and for every $k\in \{\infty,2,3\dots\}$ there are binary shifts of commutant index $k$, \cite{Pr01, Theorem 5.5}.  A binary shift has finite commutant index if and only if its bitstream is eventually periodic, \cite{BY, Theorem 5.8}.  The commutant index is a cocycle conjugacy invariant \cite{BY},\cite{Pr98, Theorem 5.5}.

\remark{Remark 2} By \cite{Pr 98a, Corollary 4.10}, all binary shifts of commutant index $2$ are cocycle conjugate, so it makes sense to ask whether, given a binary shift $\alpha$ of index $2$, the list $\{\alpha, \beta_k, k\in \Bbb N\}$ includes all binary shifts of commutant index $2$.  To see that this is not the case, note from the construction of the $\beta_k$'s in the proof of the theorem that for every $m\in \Bbb N$ the algebra generated by the first $m$ generators for $\alpha$ coincides with the algebra generated by the first $m$ generators for $\beta_k$.  Hence they have the same centers and therefore, by Theorem $3$, the same nullity sequence.  Suppose $n\in \Bbb N$ is the first positive integer for which $\nu_n = null(\Cal A_n) = 0$.  Set $c_j = 0$ for $0 \leq j \leq n-1$ and $c_{n} = 1$.  Then by \cite{Pr98} $c_0,\dots,c_n$ may be completed to form a bitstream which corresponds to a binary shift $\gamma$ of commutant index $2$. Note, however, that the corresponding Toeplitz matrix $\Cal C_n$ is the zero matrix, with nullity $n$.  Since $null(C_n) \neq null(A_n)= 0$, $\gamma$ is not in the list $\alpha, \beta_k, k\in \Bbb N$.
\endremark

\remark{Remark 3}In \cite{Pr98} it was shown that if a binary shift $\alpha$ of finite commutant index has a nullity sequence which agrees except at finitely many places with the nullity sequence of $\alpha_{infty}$, then $\alpha$ and $\alpha_{\infty}$ are cocycle conjugate.
\endremark

\subhead
Problems
\endsubhead
\roster
\item   A binary shift $\alpha$ is said to be of infinite commutant index if $\alpha^k(R)'\cap R = \Bbb CI$ for all $k$.  Is there an $\alpha$'s of infinite commutant index which is cocycle conjugate to uncountably many binary shifts?
\item   Are any two binary shifts of infinite commutant index cocycle conjugate?
\endroster

\enddocument